\documentclass[12pt,a4paper,oneside]{article}
\usepackage{amsmath,amssymb}
\usepackage{graphicx}
\usepackage{subfigure}
\usepackage{authblk}

\def\1{\mathbf{1}}
\def\<{\langle} \def\>{\rangle}

\usepackage{amsthm}

\usepackage{algpseudocode,algorithm,algorithmicx} 

\usepackage{mathtools}  

\usepackage{amsmath} 
\usepackage{amssymb} 
\usepackage{amsthm}

\usepackage{array}
\usepackage{url}
\usepackage[T1]{fontenc}
\usepackage[utf8]{inputenc}
\usepackage{color}
\usepackage{bbold}

\usepackage{subfigure}
\usepackage{caption}

\usepackage{tikz}
\usepackage{MnSymbol}  

\newcommand{\rem}[1]{}
\begin{document}
	
\title{Analyzing Raman Spectral Data without Separabiliy Assumption}

\author[1,2]{Konstantin Fackeldey}
\author[2]{Jonas R\"ohm}
\author[3]{Amir Niknejad}
\author[1,4]{Surahit Chewle}
\author[1]{Marcus Weber}
\affil[1]{Zuse Instutute Berlin, Takustra{\ss}e 7, D-14195 Berlin}
\affil[2]{Technical University Berlin, Stra{\ss}e des 17. Juni 135, D-10623 Berlin}
\affil[3]{College of Mount Saint Vincent, 6301 Riverdale Ave, New York 10471, USA}
\affil[4]{Bundesanstalt für Materialforschung und -pr{\"u}fung,  Unter den Eichen 87, D-12205 Berlin}


\maketitle

\begin{abstract}
Raman spectroscopy is a well established tool for the analysis of vibration spectra,
which then allow for the determination of individual substances in a chemical sample, or for their phase transitions.
In the Time-Resolved-Raman-Sprectroscopy the vibration spectra of a chemical sample are
recorded sequentially over a time interval, such that conclusions for intermediate products
(transients) can be drawn within a chemical process.
The observed data-matrix $M$ from a Raman spectroscopy can be regarded as a matrix product of
two unknown matrices $W$ and $H$, where the first is representing the contribution of
the spectra and the latter represents the chemical spectra.
One approach for obtaining $W$ and $H$ is the non-negative matrix factorization.
We propose a novel approach, which does not need the commonly used separability
assumption.
The performance of this approach is shown on a real world chemical example.

\end{abstract}	

\section{Introduction}

In Raman spectroscopy vibrational spectra can be detected. Analysis of those spectra provides comprehension about chemical and physical properties of molecular structures, which is important in different research areas in biology, medicine and industry \cite{Ferr03,Li14,Ku08}. Nowadays, Raman spectrometers are capable to generate spectral recordings down to the femto second time scale. Such \textit{time-resolved} Raman spectroscopy allows - besides spectral recordings of stable substances - for monitoring of events like intra molecular rearrangements and chemical reactions \cite{Sang11}. We thereby obtain measured Raman spectra as a function of time, which depicts both main characteristics of an observed process: On the one hand, each measured spectrum is a fingerprint of compounds and therefore represents the intrinsic spectra of the individual species or molecular states involved in the reaction. On the other hand, the relative contributions of the involved spectra to each measured spectrum reflect the momentary composition of the sample at the corresponding time. Through the full series of generated spectra we hence draw conclusions about the kinetics of the underlying reaction process. Consequently, the central task about time-resolved Raman data analysis is deciphering the series of measured spectra with respect to the individual component spectra and their temporal evolution.   \\   
This article is organized as follows. In Section \ref{sec: NMF}, we give an overview of NMF approaches and algorithms known so far. In particular we present the separable NMF method, which found application in the approach for spectral analysis in \cite{Lu16}. Our new  NMF approach as well as the algorithmic details of the corresponding computational method are introduced in Section \ref{NMFAlgo}. In Section \ref{sec: NRI}, we present numerical results of our novel method. On the one hand, we thereby discuss recovery results for synthetic measurement data with increasing interference of the component spectra and presence of measurement noise. On the other hand, we verify the influence of the single components of our adaptable objective function through recovery results for certain choices of weighting coefficients. 

\section{Non-Negative Matrix Factorization (NMF)}
\label{sec: NMF}

From a mathematical point of view the non-negative measurement matrix $M$, which contains the discretized time-resolved Raman spectra, can be expressed as 
\begin{align}\label{eq:MatrixFactorization}
M \  = \ WH \qquad W\in \mathbb{R}_+^{n\times r}, H \in \mathbb{R}_+^{r\times m}, 
\end{align}
where the columns of $W$ represent the component spectra and 
$H$ the  course of the relative concentrations. A factorization of $M$ into the two matrices $W$ and $H$ is from the chemical point of view interesting, the matrix $W$ gives us the substances being involved in the reaction and the matrix $H$ allows inference on the speed of the reaction. Note, that this is not possible by considering only one row or column of the matrix $M$.  
Summing up, time-resolved Raman spectral data can be modeled as the product of two non-negative matrices representing the single component spectra and the underlying reaction kinetics. 
\begin{figure}[H]
 \centering
     \includegraphics[width=0.7\textwidth]{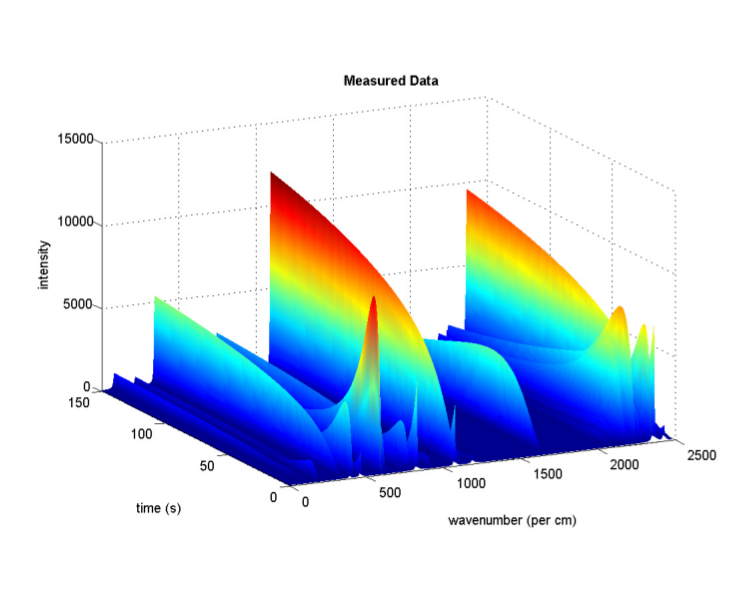}
 \caption[Interpolated visualization of the measurement data matrix $M$]{Interpolated visualization of the measurement data matrix $M$. The matrix $H$ represents the "normed" intensity which we term relative concentration. The matrix $W$ represents the wavenumber.}
  \label{Data}
\end{figure}
Recovering these factorization matrices only given the measured time-resolved spectra requires \textit{non-negative matrix factorization} (NMF). In general, NMF is an utile tool for the analysis of high-dimensional data and therefore relevant topic in present-day research in many scientific fields \cite{Gu02, Xu03, De08}. Besides detecting a compressed representation, NMF delivers insights into structure and features of the given data by extracting easily interpretable factors. 

The goal of nonegative matrix factorization (NMF) (see e.g. \cite{Gi142,De08} and the references therein) of a data matrix $M$ as input, is to solve an optimization problem 
in order to find matrices $W$ and $H$ with non-negative entries such that the product $WH$ is the best possible 
approximation of our non-negative input data matrix $M$. NMF is a linear dimension reduction technique for a non-negative data set, which means that the corresponding matrix of data points is approximated by a linear combination of the columns of matrix $W$.

\paragraph{Mathematical Background} The columns of $W$ form a basis for the column space of matrix $M$ and the columns of matrix $H$ are the weights
to approximate the data points. The NMF problem is $\cal NP$-hard \cite{Vavasis10}, due to the non-negative constraints on $W$ and $H$.
Moreover the solution of an NMF Problem is generally not unique. To see this, assume that  $W>0$, $H>0$, and that there exists a matrix $D$ such that $WD>0$ and $D^{-1}H >0$ then $M=(WD)(D^{-1}H)$ which shows that the NMF is not unique.

In the absence of the positivity constraints the problem could be solved efficiently by using methods such as truncated singular value decomposition (TSVD) \cite{H71}. One of the common approaches for solving the NMF problem is the 
{\it alternating least squares approach}  \cite{Lin2007,Berry2007}. In this approach, one of the two matrices is fixed, for example $H$ and then finds the corresponding optimal solution for $W$, which is a convex optimization problem with non-negativity constraints. Then alternate between $W$ and $H$.
If the matrix $M$ satisfies a separability condition, then we can solve the NMF problem efficiently. By definition a matrix $M$ is $r$-separable , if there exists a non-negative factorization (exact factorization) of rank $r$, where each column of $W$ is equal to a column of $M$.
Meaning that each column of $W$, being a basis for the column space of $M$, appears somewhere in the data matrix $M$ as its column. 

Geometrically, the columns of $W$ are the vertices of the convex hull of the columns of $M$. The separability condition means, that all columns of $M$ can be reconstructed by using a convex combination of $r$ columns of $W$ \cite{SeparabilityAssumption2,SeparabilityAssumption1}. This is only possible, if the columns of $M$ form a simplex which is spanned by $r$ columns of $M$. This is not necessarily the case.

\paragraph{NMF in the context of measurement data}
Given a component-wise non-negative matrix $M$ of dimension $n \times m$ and an integer $r>0$, NMF determines likewise componentwise non-negative matrices $W$ and $H$ of dimensions $n \times r$ and $r \times m$, respectively, such that $M=WH$. Generally, integer $r$ is denoted as \textit{rank} of the factorization. Assuming $M$ to represent~$m$ measurements of~$n$ non-negative variables, we interpret the NMF task as follows: We aim to identify $r$ \textit{ingredients} which allow for recovery of all $m$ measurements by composition according to respective contributions. The ingredients then are reflected by the columns of factorization matrix $W$ while the columns of $H$ contain the corresponding mixing coefficients.  \\
In practice, considering measured data and therefore allowing noise or other forms of data uncertainty generally rules out the existence of an \textit{exact} NMF in terms of $M=WH$. Thus, from now on we want to compute component-wise non-negative matrices $W$ and $H$ such that $WH$ is an approximation of $M$. 

In the context of Raman data spectral analysis, focusing on the non-negativity of involved matrices becomes reasonable through the model for time-resolved Raman spectral data of Liesen et al. \cite{Lu16}. They introduce an approach to express a series of spectral recordings of a chemical reaction (matrix $M$) as the matrix product of the component spectra (matrix $W$) and the evolution of relative concentrations of these reaction components (matrix $H$). Based on this model and synthetic spectral data, which satisfy the recently much-cited \textit{separability} assumption, the authors of  \cite{Lu16} furthermore present an algorithm to detect a factorization $WH=M$ using \textit{separable} NMF methods. 


Inspired by their results, we propose a novel method, which does \textit{not} rely on the \textit{separability} assumption, since in the context of a spectral analysis this assumption is very restrictive. The separability assumption means that the convex hull of the columns of $M$ is given by the column vectors of $W$. This is not necessarily given in real-world data. In other words, this assumption means, that the convex hull of $M$ is a simplex.  Of course, it is true that we are searching for a simplex that includes all columns vectors of $M$, but the convex hull of $M$ needs not be a simplex. Thus, we will exploit additional chemical or physical model aspects in order to find the {\em optimal} simplex including the columns of $M$ without separability assumption. 
In the center of attention of this new approach stands an adaptable objective function, taking into account only the common structural properties of the sought-for, process defining matrices $W$ and $H$. \\

\section{Solving an Optimization Problem for NMF}
\label{NMFAlgo}

In the following we pick up the concepts of both previous chapters as we introduce a new NMF approach which is specialized on analysis of time-resolved Raman spectral data. Recall from \eqref{eq:MatrixFactorization} that the thereby recovered non-negative matrices represent the component spectra of the involved species ($W$) and the reaction kinetics in terms of the evolution of relative concentrations ($H$). Our \textit{novel NMF} approach differs from the methods discussed so far as it is mainly based on minimization of an objective function which directly incorporates all known structural properties of the sought-for matrices $W$ and $H$. Furthermore, our approach is unaffected by the restrictive separability assumption. In contrast to Liesen et al. \cite{Lu16}, we hence apply our method even to non-separable measurement data. Additional flexibility and adaptability of the novel approach will be depicted in the numerical results in section \ref{sec: NRI}. Here we present the leading ideas of this approach as well as the details of the corresponding computational method. 

\subsection{Optimization Criteria for NMF}
\label{sec: MP}

In the following we propose a novel approach which is based on an objetive function 
which includes the needed structural properties of the sought-after matrices 
$W$ and $M$.

\paragraph{Claims on the matrices $W$ and $H$} In the following we assume, that the component spectra are positive, such that $W$ is a positive matrix. The componentwise non-negativity of the kinetics $H$ is also reasonable, since 
relative concentrations are in general non-negative.
Furthermore, because of representing relative concentrations, each column of $H$ is a priori supposed to sum up to 1.

For each of the $s$ chemical species the relative concentration is given by the relavtiv concentration function $h_s$:
\begin{align*}
h_s : \left[0,T \right] \rightarrow \left[ 0,1 \right] , \qquad s = 1,\dots , r. 
\end{align*}
describing the relative concentration of species $s$ at 
time $t \in \left[0,T \right]$ of the considered reaction.

Since the concentrations $h_s(t)$ are relative we have 
\[
\sum_{s=1}^r h_s(t) =1  \text{ for each  } t \in [0,T].
\]
By using $m$ time steps for discretization of the concentration functions $h_s(t)$
we obtain the column stochastic matrix 
\[
H = \begin{bmatrix} h_1(t_0)& \dots  &\dots & h_1(t_{m-1})\\
h_2(t_0)& \dots  & \dots & h_2(t_{m-1}) \\
\vdots & \dots  & \dots &\vdots \\ 
h_r(t_0)& \dots  &\dots & h_r(t_{m-1})
\end{bmatrix}.
\]

The sequential Raman-measurements can not be modelled as a ``random picking of spectra''. The temporal order of measurements  is important. Let the columns of $H$ be given by $h(t_i), i=0,...,m-1$, i.e.
\[
H= [h(t_0)|\dots \dots|h(t_{m-1})],\quad h(t_i) \in \mathbb{R}^r, i=0,...,m-1.
\]

Given the initial ``concentrations'' $h(t_{i-1})$ there is a kinetics (or some Markov process) providing the concentrations of the next time-step  $h(t_i)$. This can be modelled by assuming a transition matrix $P$ for the autonomous Markov process, if the time intervals are always constant. Thus, we claim that there exists a (row) stochastic matrix $P \in \mathbb{R}^{r\times r}$ such that
\begin{equation} \label{eq:constructP}
    (h(t_{i-1}))^T \cdot P = (h(t_{i}))^T,\qquad i=0,...,m-1.
\end{equation}

In other words, the change of the relative concentration between the time steps can be interpreted as a Markov process.  The construction of this matrix $P$ will be explained later.

Summing up the objective function in our approach has the following penalty terms

\begin{itemize}
\item[i)] $W$ is component-wise non-negative, 
\item[ii)] $H$ is component-wise non-negative,
\item[iii)] $H$ is column stochastic,
\item[iv)] $P$ is component-wise non-negative, and
\item[v)] $P$ is row stochastic.
\end{itemize}
Summing up, we arrive at the following objective function
\begin{eqnarray}\nonumber\label{eq:objectivefkt}
    \Psi & = &\alpha \left( \min\limits_{i,j} \; {W}_{ij} \right) + \beta \left( \min\limits_{i,j} \; {H}_{ij} \right)+\gamma \left( \max\limits_j  \; |\sum\limits_{i=1}^r \; {H}_{ij} -1| \right) \cr
    & & + \delta \left( \min\limits_{i,j} \; {P}_{ij} \right)+\mu \left( \max\limits_i  \; |\sum\limits_{j=1}^r \; {P}_{ij} -1| \right) 
\end{eqnarray}

It has to be mentioned here, that the constraint iv) is not necessarily valid. The matrix $P$ has to be row-stochastic, however, the entries of $P$ can be negative. A Galerkin projection of a Markov Process on the basis of microstates to a small set of macrostates can lead to negative entries in the projected matrix $P$. In the real-world example in Section \ref{sec:Ethanol}, we will show a crystallization process with a non-exponential decay of one species,
which leads to a matrix $P$ with one negative entry.     

\paragraph{Robust Perron Cluster Analysis (PCCA+)}
In the computational method of our novel NMF approach we apply the Robust Perron Cluster Analysis (PCCA+) \cite{We06} to generate an initialization of the kinetics in matrix $H$. We thus briefly introduce intention and operating principles of PCCA+ and reveal its utility for our context. \\ PCCA+ belongs to the family of algorithms for characterizing objects of similar behaviour to combine them into a certain number of clusters. In several areas of computational life science this kind of task plays a versatile role. PCCA+ arises from investigation of molecular conformation dynamics and the thereby main interest into identification of metastable conformations \cite{De05,We02}. There, metastable conformations are clusters for which the large scale geometric structure of the observed ensemble is conserved under the influence of a spatial transition operator \cite{Sch99}. Translating this approach into terms we consider a stochastic matrix $T \in \mathbb{R}^{N\times N}$ (representing the discretized version of the spatial transition operator) and we search for a non-negative matrix $Y \in \mathbb{R}^{N\times N_C}$, which column-wise contains the clusters~$y_i ,\; i=1, \dots ,N_C$, and thus satisfies three requirements: $Y$ is non-negative and row stochastic in order to meet the partition-of-unity constraint. Thirdly the vectors $y_i$  build an eigenvalue cluster near 1.0 of $T$. This means for each $i=1, \dots ,N_C$ we have  
\begin{align}
\label{eq: EV}
Ty_i \approx y_i .
\end{align} 
The main idea of PCCA+ is to generate $Y$ as a linear transformation of the matrix $X\in \mathbb{R}^{N\times N_C}$, where $X$ columnwise contains the $N_C$ first eigenvectors of $T$ with respect to eigenvalues close to $\lambda_1 = 1$. PCCA+ therefore computes a non-singular transformation matrix $\mathcal{A}\in \mathbb{R}^{N_C \times N_C}$ in order to gain the non-negative, row stochastic matrix $Y$ via
\begin{align}
\label{eq: PCCA}
Y = X\mathcal{A}.
\end{align} 
Above, in paragraph \textit{matrix properties}, we claimed that the sought-for matrix $H$ of reaction kinetics needs to be non-negative and column stochastic. Both requirements are satisfied if we consider (\ref{eq: PCCA}) and choose $H=Y^T$ as an initial guess of the kinetics. Thus, in the computational method of our novel NMF approach, the preprocessing prepares the application of PCCA+ in order to generate a promising initialization of $H$. \\ 

Investigating (\ref{eq: PCCA}) generally we may find several feasible solutions $\mathcal{A}\in \mathbb{R}^{N_C \times N_C}$ providing an appropriate matrix $Y$. PCCA+ tackles this issue by computing $\mathcal{A}$ through solving an optimization problem with respect to a certain objective function. Given that the stochastic matrix $T$ is the discretization of a transition operator (consider e.g.~molecular conformation dynamics), maximization of this objective function is equivalent to the maximization of metastability between the generated clusters. In other contexts (consider e.g.~geometrical cluster problems) the interpretation of the objective functional may be different while still meaningful. See \cite{We05,De05,We15} for exemplary applications and illustrations of PCCA+ in several research areas.

\subsection{Computational Method}
\label{Algo}
The main work stages in the computational method of our novel NMF approach are summarized in  Algorithm \ref{alg uNMF}. Note that we distinguish between the finally recovered matrices (denoted as $W_{rec}$ and $H_{rec}$) and their corresponding interim results (denoted as $\widetilde{W}$ and $\widetilde{H}$). Furthermore, we use matlab method \textit{pinv} to calculate pseudoinverses of singular or even non-square matrices. We then label the pseudoinverse of a matrix $A$ as~$A^{\dagger}$. Furthermore, with $A_+$ we denote the matrix which is constructed out of $A$ by deleting the first row and $A_-$ is the corresponding matrix constructed out of $A$ by deleting the last row. 
\begin{algorithm}
  \caption{ \quad Novel NMF for Raman Data Spectral Analysis
    \label{alg uNMF}}
  \begin{algorithmic}[1]
    \Require{data matrix $M\in \mathbb{R}^{n\times m}$ and factorization rank $r$}
    \Ensure{matrix $W_{rec}\in \mathbb{R}^{n\times r}$ of component spectra and $H_{rec}\in \mathbb{R}^{r\times m}$ of reaction kinetics such that $M \approx W_{rec}H_{rec}$} 
    \item[]
    \State{Perform SVD for primary factorization $M^T = U\Sigma V^T$ and reshape $U$ into $\mathcal{U}$.} \label{st: Pre}
    \State{Apply PCCA+ in order to initialize $\widetilde{H}=\left( \mathcal{U} \mathcal{A}\right)^T$, $\widetilde{W}  = M  \widetilde{H}^{\dagger}=M \left(\mathcal{A}^T \mathcal{U}^T \right)^{\dagger}$ and $\widetilde{P}=((\widetilde{H}_-)^T)^{\dagger}(\widetilde{H}_+)^T=({\cal U}_-{\cal A})^\dagger({\cal U}_+{\cal A})={\cal A}^{-1}\big({\cal U}_-^\dagger{\cal U}_+\big){\cal A}$.} \label{st: PCCA}
    \State{Minimize objective function with respect to transformation matrix $\mathcal{A}$.} \label{st: OF}
    \State{Reconstruct spectra $W_{rec}$ and kinetics $H_{rec}, P_{rec}$ according to the result of Step \ref{st: OF}.} \label{st: Re}
  \end{algorithmic}
\end{algorithm}

\begin{itemize}
\item \textbf{Step 1: Preprocessing} \qquad In the preprocessing we consider $M^T$. By subtraction of a reference point we transfer the columns of $M^T$ into a linear space. Afterwards we perform singular value decomposition (SVD) such that we gain $M^T = U\Sigma V^T$. In order to initialize $\widetilde{H}$ we want to apply PCCA+ to the leading $r-1$ columns of $U$. Thus we build a matrix $\mathcal{U}$, which takes the role of $X$ in (\ref{eq: PCCA}), as follows: The first column of $\mathcal{U}$ is equal to $e=\left[ 1, \dots , 1 \right]^T \in \mathbb{R}^m$, which is a requirement of PCCA+. We then stock up with columns $1, \dots , r-1$ of~$U$ until $\mathcal{U}\in \mathbb{R}^{m\times r}$. Subsequently, for efficiency reasons of PCCA+, we ensure orthogonality among the columns of $\mathcal{U}$ \cite{We06}.
\item \textbf{Step 2: Initializing $\widetilde{H}$, $\widetilde{W}$, and $\widetilde{P}$} \qquad We apply PCCA+ to $\mathcal{U}$. According to~(\ref{eq: PCCA}), we obtain a non-negative, column stochastic matrix $\widetilde{H}$ setting
\begin{align}
\label{eq: guessH}
\widetilde{H} = \left( \mathcal{U} \mathcal{A}\right)^T \; \in\mathbb{R}^{r\times m} ,
\end{align}
whereby $\mathcal{A}\in\mathbb{R}^{r\times r}$ is the computed PCCA+ transformation matrix. $\widetilde{H}$ is our initial guess of the kinetics of relative concentrations. Accordingly we gain an initialization of the component spectra $\widetilde{W}$ through the relation
\begin{align}
\label{eq: guessW}
M & = \widetilde{W} \widetilde{H} \notag \\
\Leftrightarrow \qquad \widetilde{W} & = \ M \widetilde{H}^{\dagger} =  M \left(\mathcal{A}^T \mathcal{U}^T \right)^{\dagger} \; \in\mathbb{R}^{n\times r} .
\end{align} 
In (\ref{eq:constructP}), we can see that the matrix $\widetilde{P}$ is given by 
\begin{align}
\label{eq: guessP}
 \widetilde{P} & =((\widetilde{H}_-)^T)^{\dagger}(\widetilde{H}_+)^T \notag \\
               & = {\cal A}^{-1}\big({\cal U}_-^\dagger{\cal U}_+\big){\cal A}. 
\end{align} 
Regarding (\ref{eq: guessH}), (\ref{eq: guessW}), and (\ref{eq: guessP}) we express the initial guesses of the sought-for matrices only in terms of the given and processed data ($M$, $\mathcal{U}$) and the PCCA+ transformation matrix ($\mathcal{A}$).   
\item \textbf{Step 3: Minimizing objective function}  \qquad The objective function of our novel NMF approach only incorporates structural properties of the sought-for matrices as discussed above in paragraph \textit{matrix properties}. With respect to each property we estimate a penalty value as stated in the following expressions:
\begin{equation}
  \left.
  \begin{aligned}
  \text{Penalty 1:} \qquad & \alpha \left( \min\limits_{i,j} \; \widetilde{W}_{ij} \right) \qquad \qquad \\
    \text{Penalty 2:} \qquad & \beta \left( \min\limits_{i,j} \; \widetilde{H}_{ij} \right) \qquad \qquad \\
    \text{Penalty 3:}  \qquad & \gamma \left( \max\limits_j  \; |\sum\limits_{i=1}^r \; \widetilde{H}_{ij} -1| \right) \qquad \qquad  \\   
    \text{Penalty 4:} \qquad & \delta \left( \min\limits_{i,j} \; \widetilde{P}_{ij} \right) \qquad \qquad \\
    \text{Penalty 5:}  \qquad & \mu \left( \max\limits_j  \; |\sum\limits_{j=1}^r \; \widetilde{P}_{ij} -1| \right) \qquad \qquad  \\     \end{aligned}
  \right\}
  \label{eq: PenT}
\end{equation}
In regard to non-negativity of light intensities and relative concentrations, penalties 1, 2, and 4 determine the smallest entries in matrices $\widetilde{W}$, $\widetilde{H}$, and $\widetilde{P}$. As the sum of penalty values is supposed to increase if these smallest entries appear to be negative, weighting coefficients $\alpha$, $\beta$, and $\delta$ are generally chosen negative, too. 
The requirement on $\widetilde{H}$ to be column stochastic is regarded by computing the maximal deviation of a column sum from being equal to 1.0 in penalty 3. 
Whereas, the requirement on $\widetilde{P}$ to be row stochastic is regarded by computing the maximal deviation of a column sum from being equal to 1.0 in penalty 5. 

Consider $\Psi$ to represent the sum of penalty values. As we choose the relations~(\ref{eq: guessH}) and~(\ref{eq: guessW}) for initialization, the input arguments for the objective function are the matrices $M$, $\mathcal{U}$ and $\mathcal{A}$. Since we perform optimization with respect to parameter $\mathcal{A}$, the minimization problem can be written in the form
\begin{align*}
\min\limits_{\mathcal{A}\in\mathbb{R}^{r\times r}} \; \Psi^2 .
\end{align*}
Minimizing~$\Psi^2$ hence numerically adjusts matrices $\widetilde{W}$ and~$\widetilde{H}$ according to the claimed structural properties. For computation we apply matlab method \textit{fminsearch}, which uses the simplex search method of Lagarias et al. \cite{La98}.
\item \textbf{Step 4: Recovering $W_{rec}$, $H_{rec}$, and $P_{rec}$} \qquad The minimization in Step \ref{st: OF} finally returns a transformation matrix $\mathcal{A}_{\text{opt}}$. We then recover the resulting kinetics $P_{rec}$ of relative concentrations~$H_{rec}$ and the component spectra~$W_{rec}$ according to~(\ref{eq: guessH})-(\ref{eq: guessP}) as
\begin{align*}
H_{rec} & = \left( \mathcal{U} \mathcal{A}_{\text{opt}} \right)^T = \mathcal{A}_{\text{opt}}^T \mathcal{U}^T \; \in\mathbb{R}^{r\times m} , \\
W_{rec} &= M H_{rec}^{\dagger} =  M \left(\mathcal{A}_{\text{opt}}^T \mathcal{U}^T \right)^{\dagger} \; \in\mathbb{R}^{n\times r}, \\ 
P_{rec} &= {\cal A_{\text{opt}}}^{-1}\big({\cal U}_-^\dagger{\cal U}_+\big){\cal A_{\text{opt}}} \; \in\mathbb{R}^{r\times r}. 
\end{align*} 
\end{itemize}

In regard to NMF in the context of Raman data spectral analysis, our novel approach offers two main advancements: Firstly, in contrast to the method of Liesen et al. \cite{Lu16}, our novel NMF approach is unaffected by the separability assumption. Since we only consider the general properties of the sought-for matrices without  further demands on the input data, we may apply the novel approach to the broader range of even non-separable spectral data. Secondly, note the possibility to manipulate the decicive objective function in Step \ref{st: OF} by the choice of weighting coefficients $\alpha, \beta,\gamma,\delta$ and $\mu$ or by addition of further penalty terms. This flexibility and adaptability of our method allows for example for special focus on certain data properties or even extension of the recovery objectives. We remark that the approach of optimizing $P_rec$ has already been suggested in \cite{Weber2018} and recently \eqref{eq: guessP} has been applied in \cite{GPNH_20}. \\
The next section presents some numerical experiments.

\section{Numerical Results}
\label{sec: NRI}

In this section we present the level of performance of our novel NMF approach by applying it to a sequence of artificial time-resolved Raman spectral data. After describing the reaction data generation in Section \ref{sec: RDG}, we prove that the component spectra are recovered to a high quality and that we even reach meaningful approximations of the underlying reaction kinetics. As well in Section \ref{sec: NR}, we present the effectiveness of our method in the case of increased overlap among the individual component spectra and the occurrence of measurement noise. In Section \ref{sec:Ethanol}, we present real-word data from Raman spectroscopy measured during a crystallization process of paracetamol in ethanol. We show that our method can help to identify and characterize intermediate states (and their life-times) of a chemical process. 

\subsection{Description of the Reaction Data Generation}
\label{sec: RDG}

\begin{figure}[H] 
 \centering
     \includegraphics[trim = 20mm 190mm 20mm 0mm, width=0.95\textwidth]{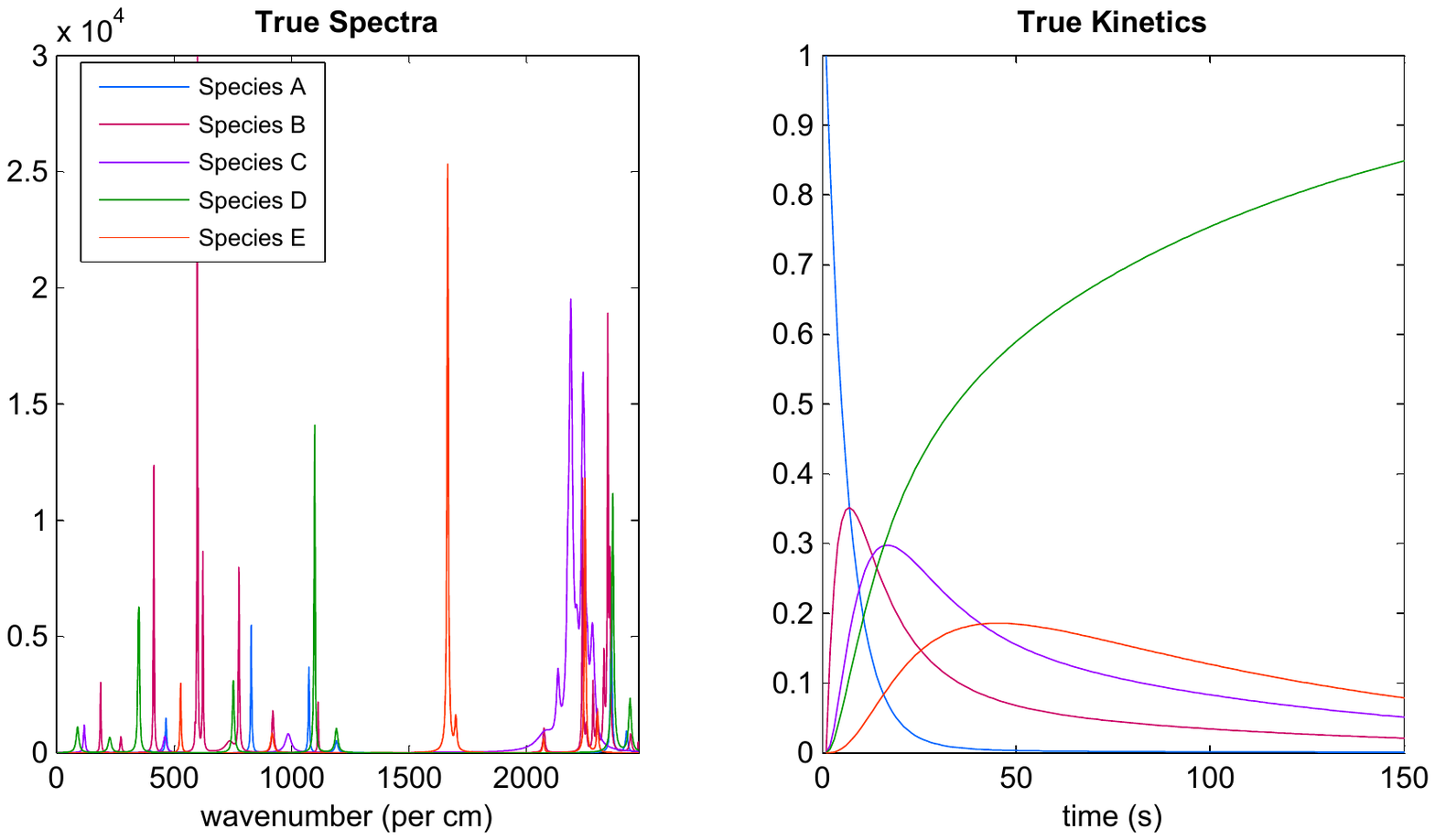} 
 \caption[Artificially generated component spectra and reaction kinetics]{Illustration of artificially generated component spectra (left) and kinetics of first-order reactions (right) including five species A to E. The assignment of color to species holds for both panels. The resulting time-resolved measurement data are displayed in Figure \ref{Data} (top).}
 \label{True}
\end{figure}

As in Section \ref{sec: NMF} for the model of time-resolved Raman spectral data, we here again follow the framework of Liesen et al. \cite{Lu16}. \\
Regarding the generation of artificial time-resolved Raman spectral data we consider a reaction scheme with five involved species A, B, C, D and E which are inter-related by first-order reactions. These first-order reactions are characterized by a rate matrix of transition coefficients as follows:   
\begin{align*}
K \;  = \;  \begin{bmatrix}
-0.53 & 0.53 & 0 & 0 & 0 \\ 0.02 & -0.66 & 0.43 & 0.21 & 0 \\ 0 & 0.25 & -0.36 & 0 & 0.11 \\ 0 & 0 & 0 & 0 & 0 \\ 0 & 0 & 0.1 & 0 & -0.1 
\end{bmatrix}
\end{align*}
The rows $i =1, \dots ,5$ of $K$ reflect the transition behaviour of the corresponding species in the course of the observed reaction. So $K_{12}$ says that 53\% of the amount of species A merge into species B per arbitrary unit of reciprocal time. The diagonal entries of $K$ represent the sum of relative loss of each species per time unit. Thus we already notice species D to be the only product of this modeled reaction as just this species exclusively absorbs rates. Here, we let species A be the only educt of the reaction and therefore denote the initial concentration vector as $h_0 \coloneqq h(t_0) =~\left[ 1, 0, 0, 0, 0 \right]^T$. With $h_0$ and rate matrix $K$ we obtain the reaction kinetics as a function of time by 
\begin{align*}
h(t)^T = \left[ h_1 (t), \dots , h_5 (t) \right] = h_0^T e^{Kt},
\end{align*}  
where $h_i (t)$ denotes the relative concentration of species $i$ at time $t$. The resulting kinetics are displayed in Figure \ref{True} (right). We gain the corresponding matrix~$H$ of kinetics by discretization of $h(t)$ at equidistant time steps $t_0, \dots , t_{m-1}$ such that $H=~\left[ h(t_0), \dots , h(t_{m-1}) \right]$. \\
The single component spectra are built up as arbitrary sums of Lorentzians, which we illustrate in Figure \ref{True} (left). The five columns of matrix $W$ accordingly contain the discretized \textit{intensity}-by-\textit{wavenumber} signals.
 The spectral overlap among the single component spectra is adjustable. This means we may increase the level of spectral interference by moving all base points $x_0$ of the generated Lorentzians towards certain focal points. The level of spectral interference decides the level of separability of the measurement data. While the results in~\cite{Lu16} are based on near-separability because of low spectral interference, we prove the effectiveness of our method even in the case of high interference among the component spectra. \\
The resulting measurement data matrix $M$ is obtained as the product of matrix $W$ of component spectra and matrix $H$ of the underlying reaction kinetics as $M=WH$. See Figure \ref{Data} (top) for an interpolated visualizatoin of $M$.

\subsection{Recovery Results}
\label{sec: NR}

Considering the measurement data according to the artificial reaction scheme as introduced in the previous Section \ref{sec: RDG}, our goal is now to recover the single component spectra as well as the reaction kinetics only given matrix~$M$. In other words, we compute matrices~$W_{rec}$ and~$H_{rec}$ by applying our novel NMF approach to $M$. We thereby are especially interested into the reconstruction of the true component spectra $W$ in order to provide a powerful tool for compound identification in real-life Raman spectral analysis. Recall that the objective function in our approach is based on adding up the penalty terms in (\ref{eq: PenT}), which represent the structural properties of the sought-for matrices and which are weighted by choice of the coefficients $\alpha , \beta$ and $\gamma$. In this section we present results of our method for the predefinitions 
\begin{align}
\label{eq: set1}
\alpha  = -0.0001, \; \; \beta  = -1 \;  \text{ and } \; \gamma = 1.
\end{align}
Recall additionally that we applied singular value decomposition in the preprocessing of our computational method. That is why the order of species in the recovered matrices $W_{rec}$ and $H_{rec}$ may be permuted in comparison to the order in the exact matrices $W$ and $H$. For comparative visualization of our recovery results we thus compute the correlation coefficients between the columns~(~$\sim$~species) of $W_{rec}$ and $W$ and associate the spectra as well as the reaction kinetics according to the maximal correlation values. \\

\begin{figure}[H]
 \centering
      \includegraphics[trim = 20mm 25mm 15mm 25mm, clip, width=\textwidth,height=0.85\textheight]{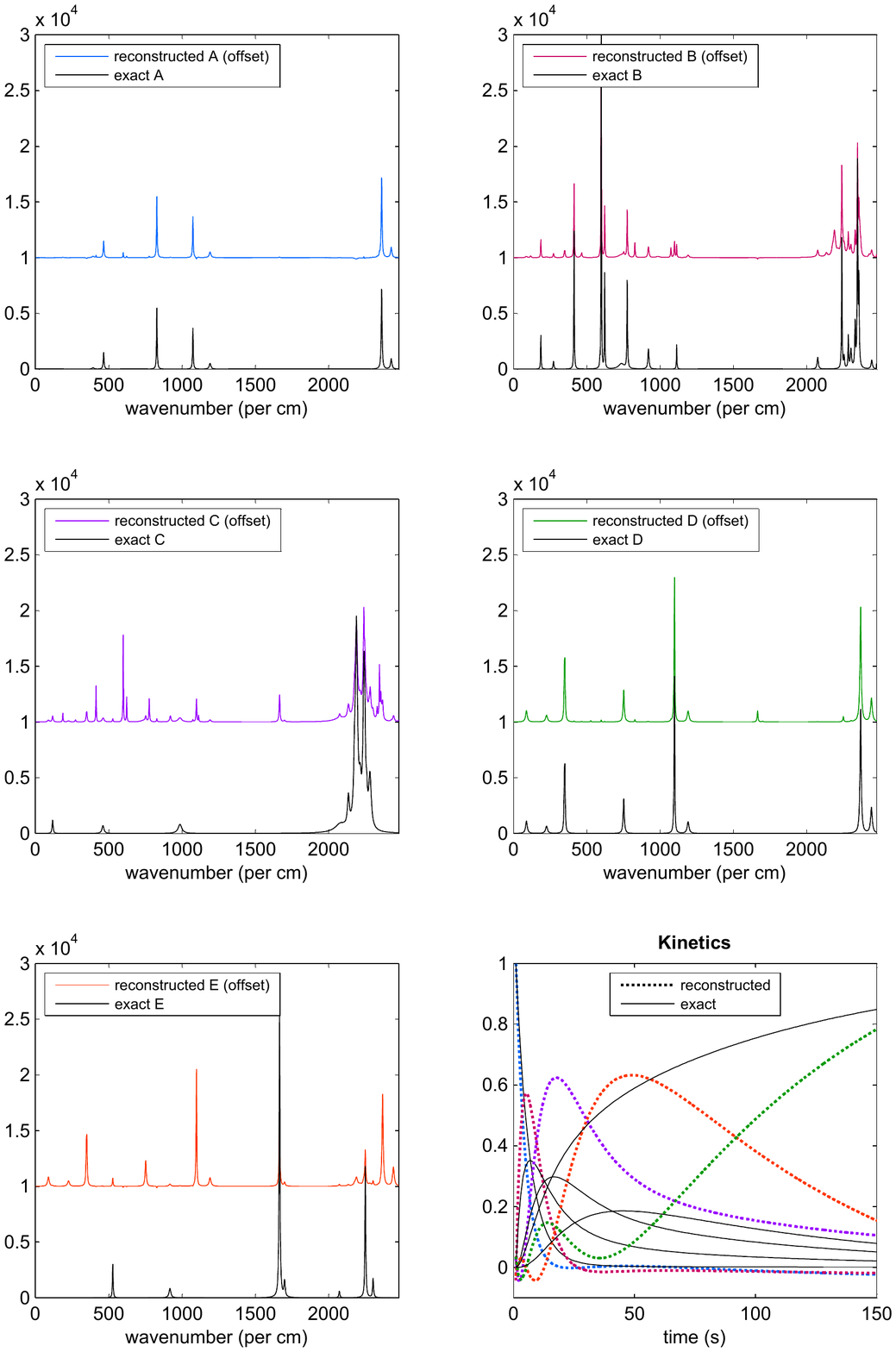}
 \caption[Recovery results for noiseless measurement data]{Reconstructed component spectra of the single species and reaction kinetics (bottom right) for noiseless Raman data. The spectra of compounds A, B and D are recovered nearly exactly. Inaccuracies in the lower wavenumber regions occur for compounds C and E. Furthermore, our computed kinetics reflect the rough behaviour of the real kinetics.}  
  \label{oSoN}
\end{figure}

Exemplary recovery results of our novel method for the noiseless case with low spectral interference are displayed in Figure \ref{oSoN}. Especially the recovery of components A, B and D is nearly exactly: The coordinates as well as the heights of peaks can hardly be distinguished visually from the original data. In the bottom right panel we also present the recovery result for the matrix $H$ of reaction kinetics.
\begin{figure}[H]
 \centering
     \includegraphics[trim = 10mm 15mm 10mm 15mm, clip, width=\textwidth]{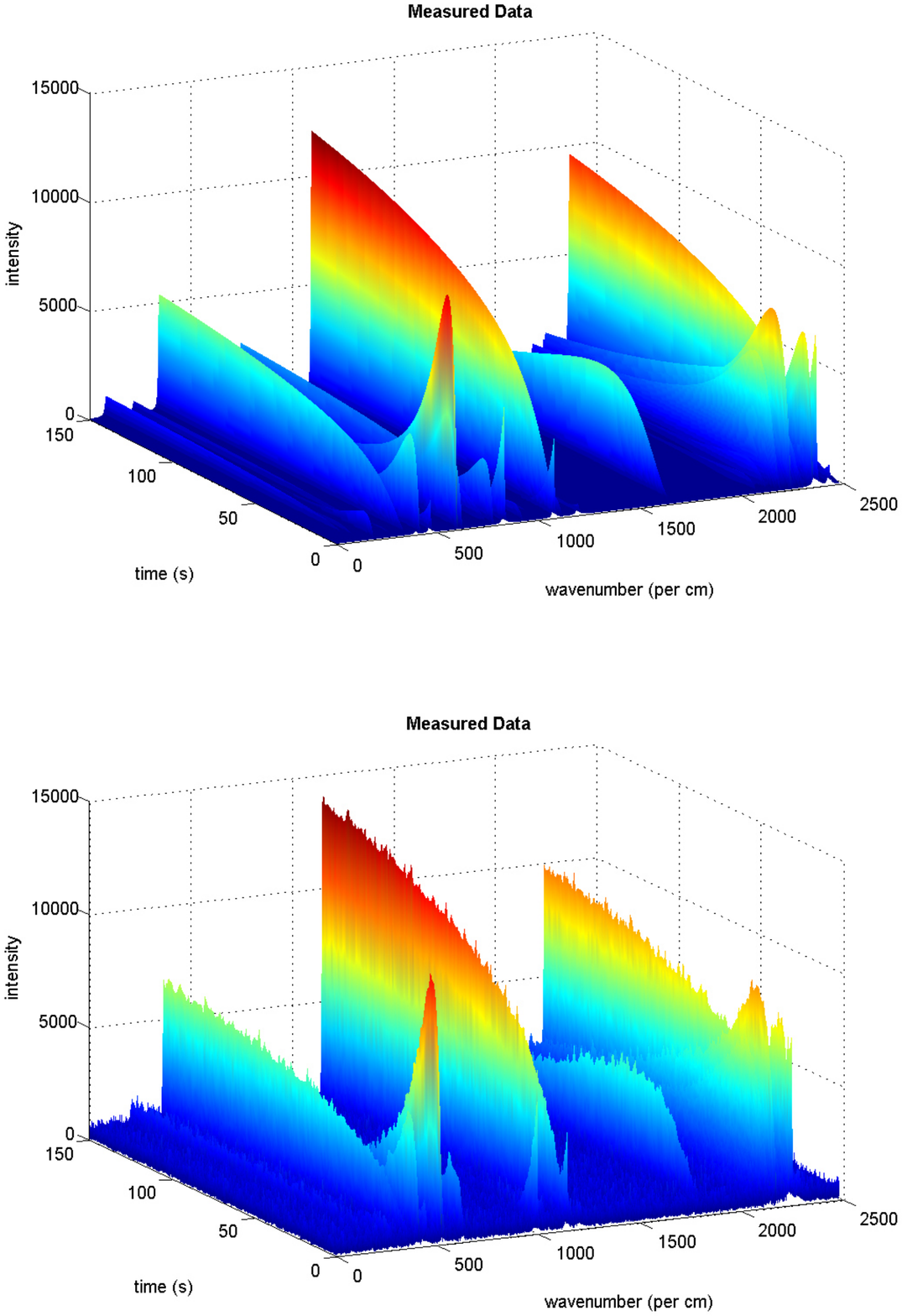}
 \caption[Interpolated visualization of the measurement data matrix $M$]{Interpolated visualization of the measurement data matrix $M$: On top, the case of well separation of the component spectra and no measurement noise. Below, a variant of increased spectral interference and noise contamination.}
  \label{Data}
\end{figure}

As in all upcoming illustrations of the reconstructed kinetics the dotted lines are assigned to their species through the corresponding color in the spectral panels. For comparison, the \textit{exact} kinetics (black lines) represent the kinetics from Figure~\ref{True}~(right). Indeed our reconstructed kinetics in Figure \ref{oSoN} reflect the general trends of the exact kinetics as in particular species A is recognized to be the only educt and species D to be the exclusive product of the generated reaction scheme.  \\

\begin{figure}[H]  
 \centering
     \includegraphics[trim = 0mm 185mm 0mm 0mm, width=0.95\textwidth]{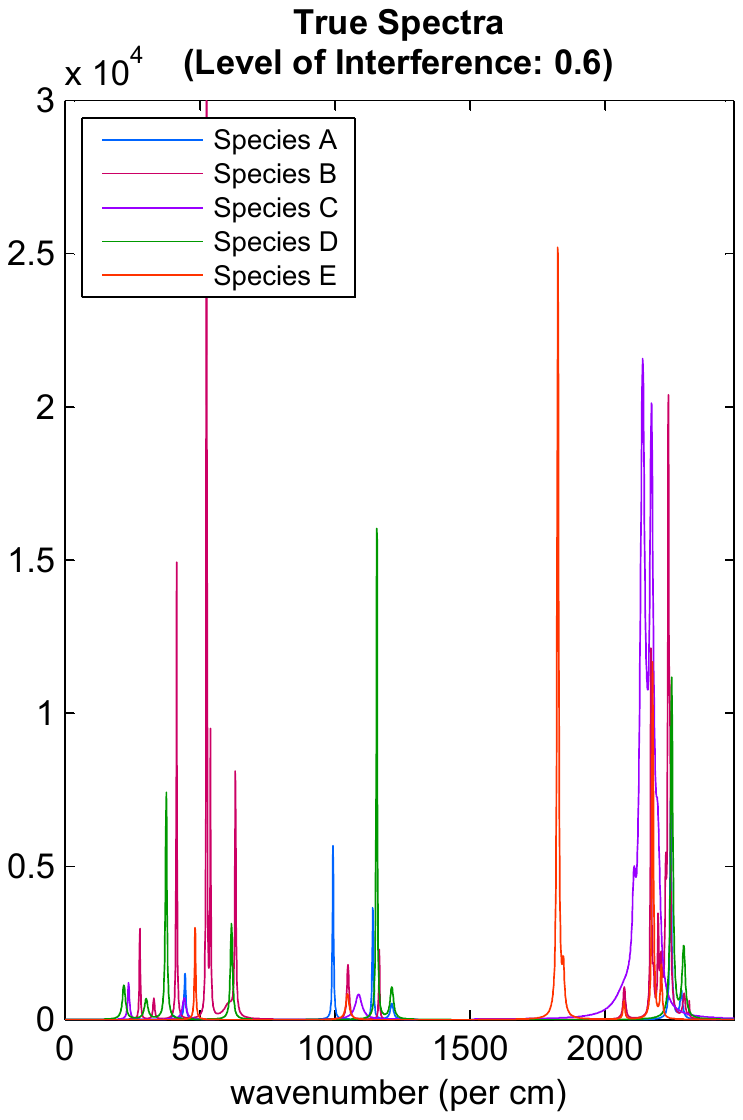} 
 \caption[True spectra for increased spectral interference]{Component spectra for modest spectral interference. In comparison to the spectra in Figure \ref{True} (left), notice how the base points of the Lorentzians have been moved closer to each other.}
 \label{SpecInt}
\end{figure}

\begin{figure}[H] 
 \centering
      \includegraphics[trim = 20mm 25mm 15mm 20mm, clip, width=\textwidth,height=0.9\textheight]{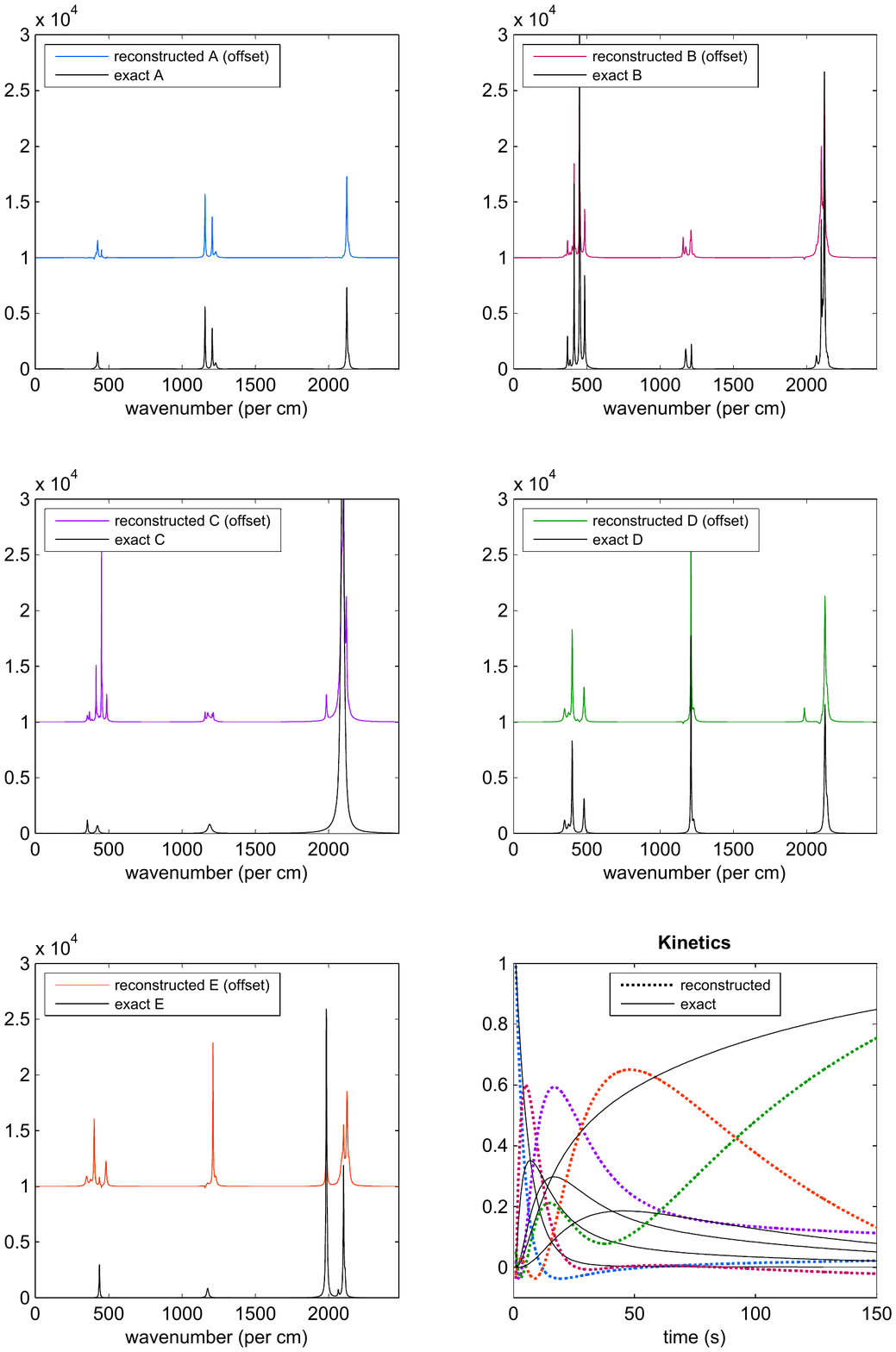}
 \caption[Recovery results for noiseless, interference-rich measurement data]{Reconstructed component spectra of the single species and reaction kinetics (bottom right) for the case of high spectral interference. Note the improvements in the recovery of species C and E in comparison to Figure~\ref{oSoN}. In addition, the educt and the product of the reaction are clearly recognizable in the recovery of reaction kinetics.} 
  \label{mSoN}
\end{figure}

As the first extension of the data setting we now investigate the effectiveness of our method in the case of increased spectral interference. As mentioned in Section \ref{sec: RDG}, we generate increased spectral interference among the component spectra in $W$ by moving the base points $x_0$ in all species towards three focal points. We then obtain component spectra as displayed in Figure \ref{SpecInt}.   
In Figure \ref{mSoN} we present the results of our novel approach being applied to very interference-rich measurement data. Besides the remaining high quality in the recovery of components A, B and D the reconstruction of species C and E apparently improved compared to the results in Figure \ref{oSoN}. In this interference-rich case our method computes the coordinates of the peaks in all component spectra quite satisfactorily. Concerning the recovery of the reaction kinetics, displayed in the bottom right panel, we again precisely identify the educt and the product of the reaction.

\begin{figure}[H] 
 \centering
      \includegraphics[trim = 20mm 25mm 15mm 25mm, clip, width=\textwidth,height=0.9\textheight]{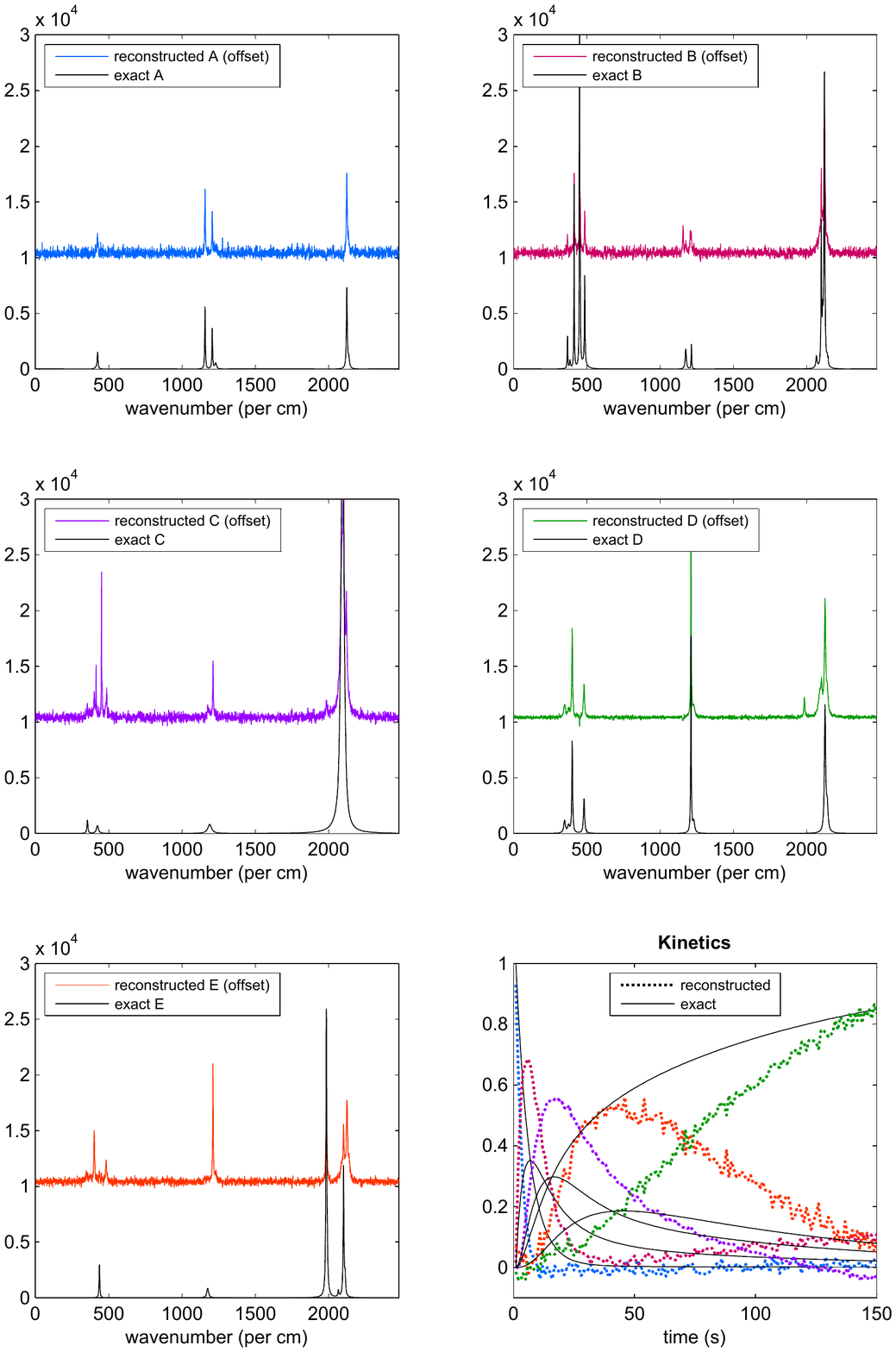}
 \caption[Recovery results for noisy, interference-rich measurement data]{Reconstructed 
 component spectra of the single species and reaction 
 kinetics (bottom right) for interference-rich and 
 noisy measurement data. The spectral recoveries still show a reasonable agreement with the true spectra. The main traits of the reaction kinetics are recognizable as well.}  
  \label{mSmN}
\end{figure}

As the second extension of our data setting we regard 
the recovery results of our routine additionally 
considering contamination of measurement noise. 
In any practical setting Raman spectral analysis needs
to deal with this issue since, for instance, signal shot noise or background noise appear in any real 
experimental data. Here we assume the noise from all 
different sources to be adequately represented by 
additive Gaussian white noise, which disturbs the 
measurement matrix $M$ according to
\begin{align*}
\tilde{M} = M + \delta \; \text{abs} \left( N \right) .
\end{align*}  
The entries of $N$ thereby are generated by the normal distribution $\mathcal{N} (0,1)$ and $\delta=0.5$ is 
the relative noise level. See Figure \ref{Data} (bottom) for an interpolated visualization of the 
interference-rich and noisy measurement matrix $\tilde{M}$. Applying our novel NMF approach with the 
predefinitions in (\ref{eq: set1}) to $\tilde{M}$, the illustrations of results in Figure \ref{mSmN} prove that the component spectra still show a reasonable agreement with the exact spectra. Furthermore, the main traits of the true reaction kinetics are recognizable in the recovered kinetics as well.  

\subsection{Example: Paracetamol in Ethanol} \label{sec:Ethanol}
\begin{figure}[H]  
 \centering
     \includegraphics[width=0.95\textwidth]{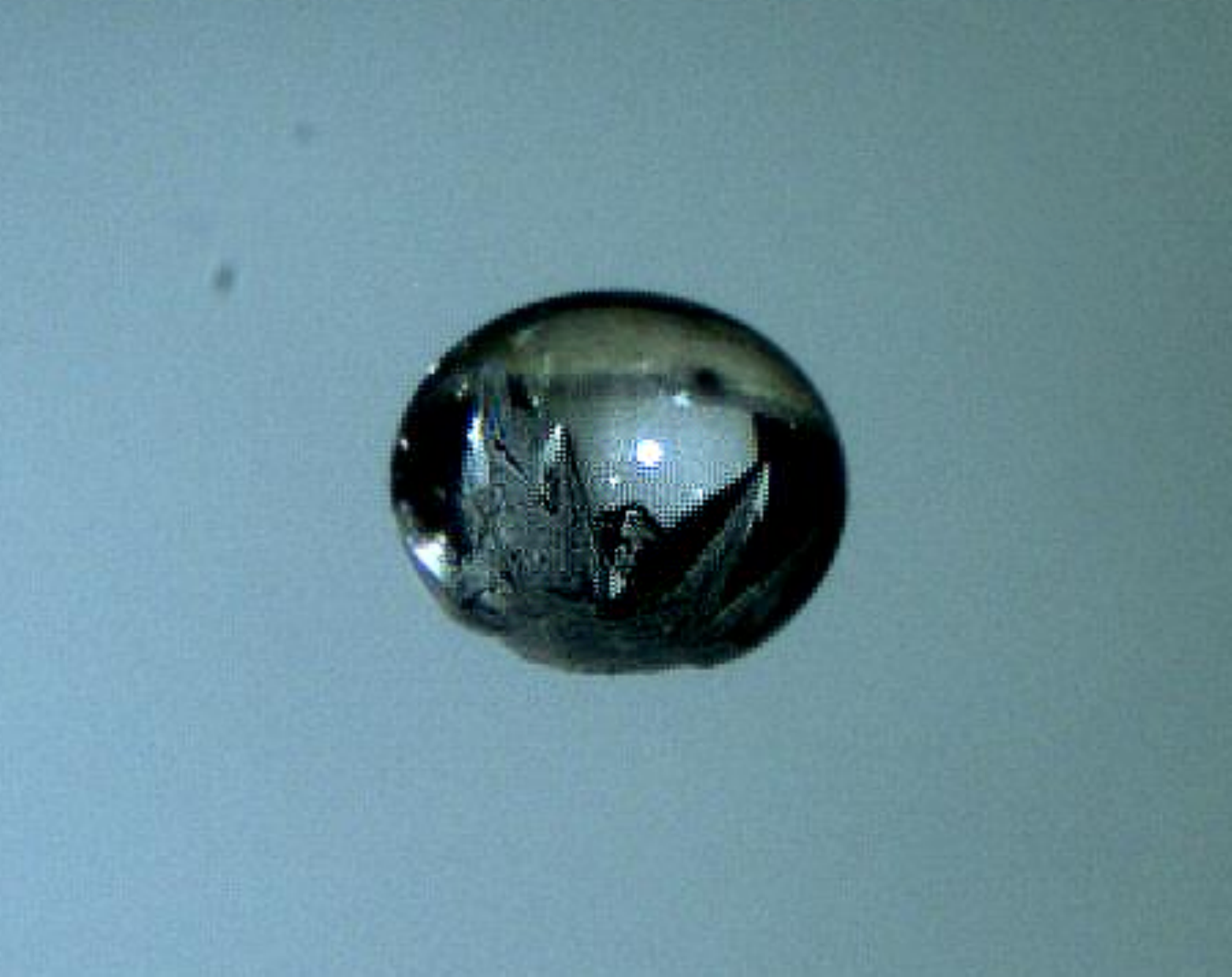} 
 \caption{Paracetamol polymorph type I crystallizing in acoustically levitated droplet of its supersaturated solution in ethanol.}
 \label{ac_img}
\end{figure}

\begin{figure}[H]  
 \centering
     \includegraphics[width=0.95\textwidth]{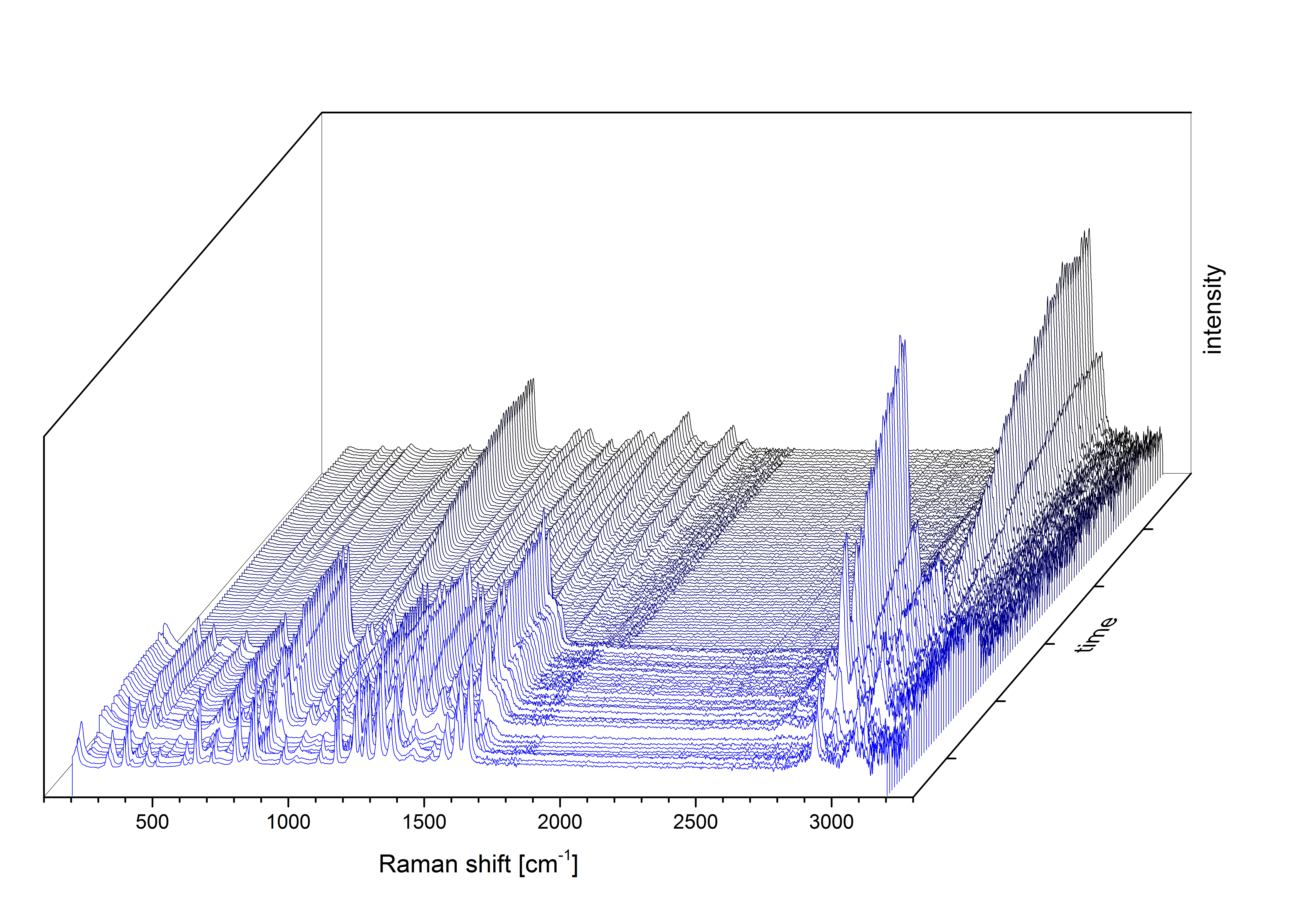} 
 \caption{In real-world applications, sequential measurements of Raman spectra lead to input data for NMF. The intensity of different wavenumbers is measured at different timesteps.}
 \label{SpecInt1}
\end{figure}

We took experimental time-resolved Raman spectroscopy data of paracetamol as an example to demonstrate application and usability of our NMF algorithm. Paracetamol crystallizes in two polymorphs, and these polymorphs can have difference in the processing of the drug in its final tablet formulation. The bioavailability of the drug can also be different according to a particular polymorph \cite{bauer2001ritonavir}. Control over crystallization is required in an attempt to manufacture a desired polymorph, for which crystallization is studied in an empirical manner with different solvents, cooling rate, etc. The effects of the solvents on crystallization of small drug molecules, paracetamol are of paramount importance. Different solvent choices yield different polymorphs of paracetamol \cite{Hilfiker2006}. Crystallization studies from liquid solutions were performed in a custom-made acoustic levitator \cite{Schlegel2012}. The acoustic levitator allows executing contact-free crystallization studies and \emph{in situ} measurements. The droplet of the solution can be fixed in a stable and undisturbed position by means of an ultrasonic field. The environment around the sample can be controlled regarding the surface, temperature, and humidity by passing a cool/hot stream of nitrogen. During the experiment the solvent evaporates and leads to a gradual increase of the concentration of the droplet which finally crystallizes (Fig.~\ref{ac_img}). Time-resolved  Raman spectroscopy is performed with the resolution of 3 seconds during this crystallization process. Various pathways from solution phase of the drug molecules to final crystallized phase have been suggested. An intermediate metastable polyamorphic state has been reported wherein the paracetamol molecules existing in transient disorganised cluster undergoes ordering to fetch final crystal structure of high order \cite{NguyenThi2015}. With our method, we were able to not only understand the kinetics of the intermediate phase, but were also able to calculate the spectra of the intermediate state. This data is crucial in understanding and thus controlling the crystallization of a drug substance. The measurements are shown in Fig.~\ref{SpecInt1}.  

\begin{figure}[H]  
 \centering
     \includegraphics[width=0.95\textwidth]{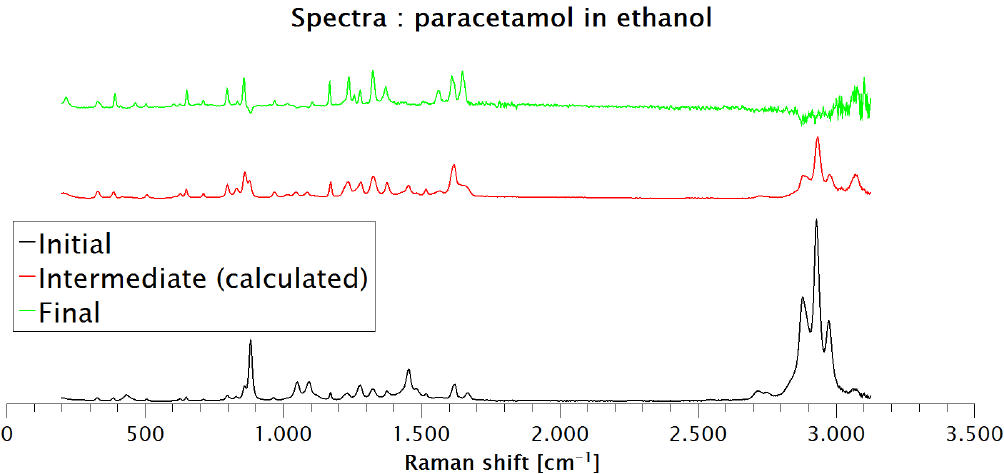} 
 \caption{During the crystallization, solvated paracetamol (black spectrum) passes through an intermediate amorphous state (red spectrum) which then immediately turns into a crystal structure (green spectrum). The three component spectra of this process are extracted by using NMF.}
 \label{SpecInt2}
\end{figure}

\begin{figure}[H]  
 \centering
     \includegraphics[width=0.95\textwidth]{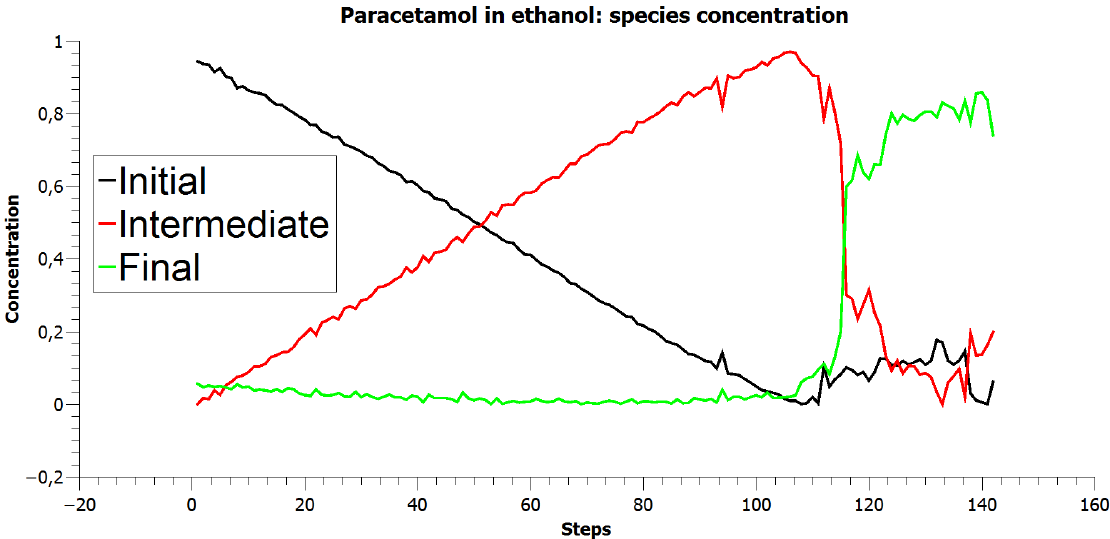} 
 \caption{Using NMF, the three ($r=3$) different compenent spectra show up during the course of time with different relative weights (concentrations). The red curve indicates initial moieties, the red curve denotes intermediate moieties, and green curve is used to indicate final crystallized polymorph. The matrix $H_{rec}$ includes the kinetics of the crystallization process.}
 \label{SpecInt3}
\end{figure}

The following settings are used for the optimization function: $\alpha=0.00001, \beta=100, \gamma=100, \delta=1, \mu=1$. With these settings it is focused on feasible concentrations. This means, we focus on providing a matrix $H_{rec}$ with non-negative entries and rowsum $1$, such that Fig.~\ref{SpecInt3} shows mathematically feasible concentration curves.  $\alpha$ is set to a very low value, because the intensities of the spectra are orders of magnitude higher than the entries in $H_{rec}$ or $P_{rec}$. After using the optimization approach Alg.~\ref{alg uNMF}, especially the matrices $H_{rec}$ and $W_{rec}$ are important experimental findings. They show the spectra of intermediate steps  and of the final crystal form of paracetamol (Fig.~\ref{SpecInt2}) and they show the kinetics of the crystallization process (Fig.~\ref{SpecInt3}). The matrix $P_{rec}$ is:
\begin{equation*}
P_{rec}=
\begin{pmatrix}  
   1.00 & 0.00 & 0.00 \\
   0.02 &  0.98 & 0.00 \\
   -0.01 &  0.02 &  0.99
\end{pmatrix}.
\end{equation*}
This matrix represents the approximated Galerkin projection (3 states) of a transition process in a continuous space (micorscopic 3D arrangement of the atoms in the droplet). The third row of $P_{rec}$ represents the initial state. The second row is the intermediate state. There is a zero probability for going back from this state to the initial state. The first row represents the stable final crystal. The upper right part of $P_{rec}$ is zero. This is because the crystallization process is directed. Fig.~\ref{SpecInt3} shows a decay of the initial state which is nearly linear. In reaction kinetics we usually expect exponential decay. The matrix is just the optimal fit to a presumed kinetics according to the chosen objective function. 
Depending on the optimization criterion, one can obtain different results from NMF of the given raw Raman spectroscopy data. These results can be checked using a cross-validation method to confirm the mathematical interpretation of the chemical process. We compared the results of NMF with simultaneous time-lapse photography of the droplet, the first of its kind to be used as a watchdog for comparing results obtained from NMF that correspond to the experimental results. Besides comparing time-step of phase change point observed in concentration curves with the experimental time-steps, another factor that validates the results are the peaks reported for metastable intermediate amorphous state closely matches with our calculated spectra. The peaks in red curve, for measured intermediate state, 1236 cm\textsuperscript{-1},1326 cm\textsuperscript{-1},1618 cm\textsuperscript{-1} to refer to few of many, match with calculated peaks at 1235 cm\textsuperscript{-1}, 1327 cm\textsuperscript{-1},1619 cm\textsuperscript{-1} \cite{NguyenThi2015}. Naturally, the peaks for final moieties can also be verified and are in accordance with reported experimental data. Structural changes, which are predicted with NMF are verified on the basis of this recording.  

\section{Conclusion}
Summarizing, our novel NMF approach returns remarkable and robust results in the recovery of component spectra and reaction kinetics while the method is mainly based on the general structural properties of the sought-for matrices. The recovery results of our approach even indicate that the quality of the recovered  component spectra improves as the spectral overlap among the component spectra increases. Our novel approach can therefore be considered as a complement to the method of Liesen et al.~\cite{Lu16} since the success of their method especially depends on low spectral interference (near-separability of $M$).

\bibliographystyle{IEEEtran}
\bibliography{Literature}

\end{document}